\documentclass[12pt]{article}
\usepackage{graphicx}
\usepackage[ascii]{inputenc}
\usepackage[T1]{fontenc}
 \usepackage[english]{babel}
\usepackage{amsmath,amssymb,amsfonts,textcomp,epsfig}
\usepackage{hyperref}
\usepackage{color}
\usepackage{calc}
 \usepackage[all]{xy}
\def\Im{ \textrm{Im}\,}
\def\R{\mathbb{R}}
\def\C{\mathbb C}
\def\Z{\mathbb Z}

\renewcommand\Re{\operatorname{Re}}
\renewcommand\Im{\operatorname{Im}}
\title{On the number of limit cycles which appear by perturbation of
 two-saddle cycles of planar vector fields}
 \author{Lubomir Gavrilov \\
 \normalsize \it Institut de Math\'{e}matiques de Toulouse, UMR 5219\\
 \it Universit\'{e}  de Toulouse,  31062 Toulouse,  France  }
\usepackage{pdfsync}
\begin{document}
\maketitle
\newtheorem{definition}{Definition}
\newtheorem{remark}{Remark}
\newtheorem{theorem}{Theorem} 
\newtheorem{lemma}{Lemma}
\newtheorem{proposition}{Proposition}
\newtheorem{corollary}{Corollary}
\vspace{5mm} \noindent 2000 MSC scheme numbers: 34C07, 37G15, 70K05
\begin{abstract}
We prove that the number of limit cycles, which bifurcate
from a two-saddle loop of an analytic plane vector field $X_0$, under an arbitrary finite-parameter analytic deformation $X_\lambda$, $\lambda \in(\mathbb{R}^N,0) $, is uniformly bounded with respect to $\lambda$.
\end{abstract}
\newpage

\section{Introduction}
Consider a finite-parameter analytic family of analytic plane vector fields
\begin{equation}
\label{xl}
X_\lambda= P(x,y,\lambda) \frac{\partial}{\partial x}+ Q(x,y,\lambda) \frac{\partial}{\partial y},\quad \lambda \in\mathbb{R}^N 
\end{equation}
such that $X_0$ has a limit periodic set $\Gamma$.  The
cyclicity  of $\Gamma$ is, roughly speaking, the
maximal number of limit cycles of $X_\lambda$ which tend to $\Gamma$ as
$\lambda \rightarrow 0$. The Roussarie's finite cyclicity conjecture claims that 
\emph{every limit periodic set occurring in an analytic finite-parameter family of plane analytic vector fields, has a finite cyclicity} \cite{rous98}. If true, the conjecture would imply the finitness of the maximal number $H(n)$ of the limit cycles, which a plane polynomial vector field of degree $n$ can have. 
Therefore it plays a fundamental  role in all questions related to the second part of the 16th Hilbert problem and its ramifications. 

\begin{figure}[htbp]
\begin{center}
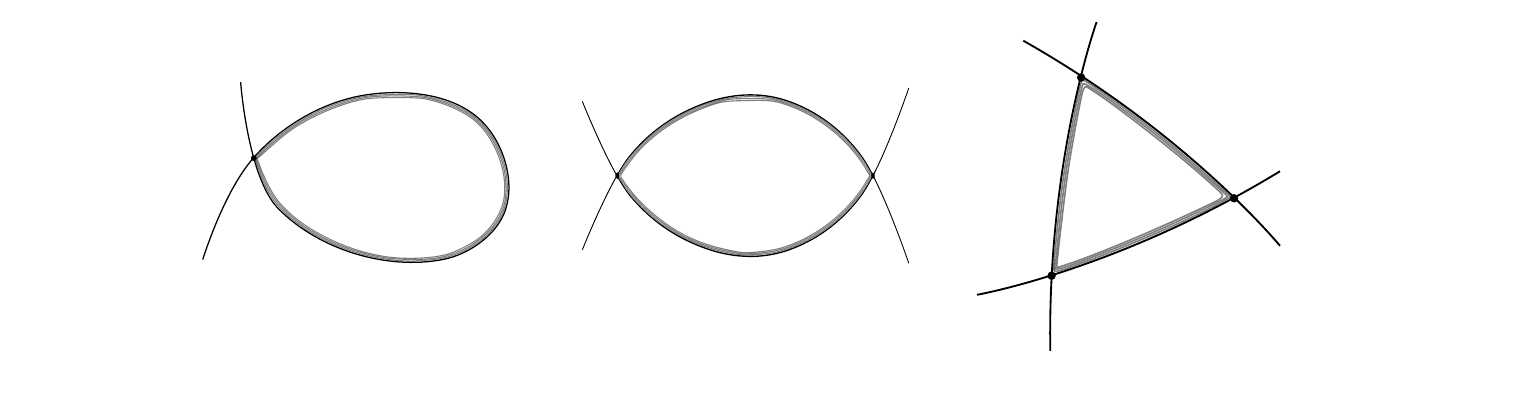
\caption{One, two and three-saddle cycles.}
\label{dessin9}
\end{center}
\end{figure}

Recall that a polycycle
of the vector field $X_0$ is a topological polygon composed of separatrices and singular points. 
A $k$-saddle cycle of $X_0$ (or a hyperbolic $k$-graphic) denoted $\Gamma_k$,
is a polycycle composed of $k$ distinct saddle-type singular points
$p_1,p_2,\dots,p_k$, $p_{k+1} =p_1$ and separatrices (heteroclinic  orbits)
connecting $p_i$ to $p_{i+1}$ as on fig. \ref{dessin9}.
$k$-saddle cycles, period orbits and weak foci or centers are the simplest limit periodic sets.
The finite cyclicity of period orbits and weak foci is well known and follows from the Gabrielov theorem \cite[p.68]{rous98}.
The finite cyclicity of one-saddle loops is due to Roussarie \cite{rous86,rous89}.

\emph{The purpose of the present paper is to prove the finite cyclicity of  a two-saddle cycle, under finite-parameter analytic deformation},  see Theorem \ref{mainth}.

Several special cases of this result were earlier proved, under different genericity assumptionss either on $X_0$ or on the family $X_\lambda$, by Cherkas, Mourtada, El Morsalani, Dumortier, Roussarie, Rousseau, Jebrane, {\.Z}o{\l}adek,Li, Caubergh, Luca and others  \cite{cher68,liro04,duro06,cdr07,ldcr09}, see also \cite[section 5.4.1]{rous98} for survey of the results and references up to 1996. The finite cyclicity of  a k-saddle cycle (any k), under finite-parameter analytic deformation was recently announced by Mourtada \cite{mour09}.

For \emph{generic families} of vector fields the analitycity can be relaxed. As it is shown by Ilyashenko and Yakovenko \cite{ilya95}, Kaloshin \cite{kalo03}
\emph{any nontrivial elementary polycycle occurring in a generic $k$ - parameter family of $C^\infty$ vector fields has finite cyclicity}. 

In contrast to the aforementioned papers we shall not use the asymptotic expansions of the corresponding Dulac maps. Instead of this, we   evaluate the number of the  limit cycles near $\Gamma_2$ in a complex domain, by making use of a suitable version of the argument principle.  This approach was initiated by the author in \cite{gavr11}, where we studied cyclicity  of Hamiltonian two-loops. As it is well known, the limit cycles of planar systems close to Hamiltonian
are closely related to the zeros of associated Abelian integrals depending on a parameter (the so called weakened 16th Hilbert problem \cite[Arnold, p.313]{arno88}). Zeros of complete elliptic integrals were successfully studied by topological arguments in a complex domain (the argument principle) after the pioneering work by G. S. Petrov \cite{petr86,petr88b}, see also {\.Z}o{\l}adek \cite[section 6]{zola06} for a description of the method. It has been used in a more general context in several papers, e.g. \cite{bmn09,gny10}, and in \cite{gavr11} the idea has been used to replace Abelian integrals by the true Poincar\'e return map.

In the present paper we shall find a relation between the fixed points of the Poincar\'e first return map and the fixed points of holonomies of the  separatrices of the saddle points, which correspond to complex limit cycles. 
To count such fixed points is a question on the zeros of families of \emph{analytic functions} which is easily solved. The main technical tool is Lemma 
\ref{mainlemma} in which we show that the connected components of the zero locus of the imaginary part of a Dulac map are smooth semianalytic curves. This allows to estimate the variation of the argument of the displacement map along the border of an appropriate complex domain, and finally to apply the argument principle in order to evaluate its zeros in the domain. 

Note that previously the relation between the monodromy and the Dulac map was used by Roussarie to compute the Bautin ideal associated to the Dulac map \cite{rous98a}. This combined with \cite{rous86, rous89} also implies the finite cyclicity of one-saddle cycles.

The paper is organized as follows. In section \ref{dulacsection} we provide the necessary technical background, and prove  the main technical Lemma  \ref{mainlemma}. 
In section \ref{onesaddle}, we give a new self-contained proof of Roussarie's theorem about the finite cyclicity of one-saddle cycles. 
The origin of our method is then explained 
in section \ref{petrovtrick} where we give a brief account of a local version of the so called "Petrov trick".
The same method is  then easily adapted in section \ref{mainsection} to show that the cyclicity of $\Gamma_2$ is finite.

\section{The Dulac map}
\label{dulacsection}

\begin{figure}
\begin{center}
\resizebox{12cm}{!}{
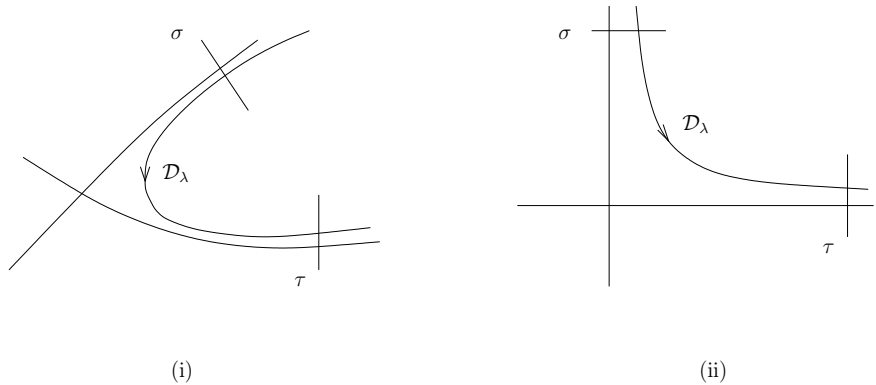}
\end{center} \caption{The Dulac map}\label{fig1a}
\end{figure}

Consider an analytic family of plane real analytic foliations $\mathcal{F}_\lambda$, $\lambda \in \mathbb{R}^N$, having a non-degenerate isolated saddle point. An appropriate translation analytically depending on $\lambda$ will place the saddle point at the origin.
The foliation $\mathcal{F}_\lambda$
has two  analytic separatrices, transversally intersecting at the saddle point, and  depending
analytically on $\lambda$ \cite{brbo56,mamo80}. Therefore a further 
 real bi-analitic change of the variables $x,y$,analytically depending on $\lambda$, 
 will identify them to the axes $\{x=0\}$
and  $\{y=0\}$ as on fig.\ref{fig1a}, so 
\begin{equation}
\label{foliation}
\mathcal{F}_\lambda : 
 x(1+...) dy +  \alpha(\lambda) y(1+
...)dx, \; \alpha(0)>0 .
\end{equation}
where the dots replace higher order terms in $x,y$ with coefficients depending on $\lambda$. 
The number $\alpha(\lambda)$ is the hyperbolic ratio of the saddle point. 
From now on we shall suppose that the foliation  (\ref{foliation}) is analytic and depends analytically in $\lambda$ in a neighborhood of the origin in $\R^2 \times \R^N$.

For $c_1,c_2\in \R$ sufficiently small, let
  $\sigma\subset \{y = c_1\}$, $\tau\subset \{x=c_2\}$  be  open complex discs centered at $(0,c_1)$ and $(c_2,0)$, parameterized by $x$ and $y$ respectively.
The (real) Dulac map is the germ of analytic map at $x=0$
$$
\mathcal{D}_\lambda : \sigma \cap \mathbb{R}^+\rightarrow \tau\cap \mathbb{R}^+, \quad \mathcal{D}_\lambda(0)=0
$$
defined as follows: if $x \in \sigma\cap \mathbb{R^+_*}$ then
$\mathcal{D}_\lambda(x)\in \tau\cap \mathbb{R^+_*}$ is the intersection with $\tau\cap \mathbb{R^+_*}$
of the orbit $\gamma_\lambda(x)$ of (\ref{foliation}), passing through $x$, see  \ref{fig1a} (ii). This geometric definition of $\mathcal{D}_\lambda$ allows to control to a
certain extent its analytic continuation in a complex domain.
\subsection{Analytic continuation}
\label{scontinuation}
 The Dulac map  allows an analytic
continuation on some open subset of the universal covering  space $\sigma_\bullet$ of
$\sigma\setminus \{0\}$, depending on $\lambda$.
Let us parameterize    $\sigma_\bullet$ by polar
coordinates $\rho > 0, \varphi \in \R$, $z= \rho \exp{i \varphi}$.
The following result is well known (e.g.  \cite[Appendix A]{gavr11})
\begin{theorem}
\label{continuation} There exists $\varepsilon_0>0$ and a continuous function
\begin{equation*}
\begin{array}{rcl}
\rho: \R& \rightarrow &\R^+_*\\
     \varphi &\mapsto & \rho(\varphi)
\end{array}
\end{equation*}
such that the Dulac map allows an analytic continuation in the domain
\begin{equation}\label{domain}
\{(\lambda, \rho, \varphi)\in \C^N\times\sigma_\bullet: |\lambda|<
\varepsilon_0, 0< \rho < \rho(\varphi)\}
\end{equation}
\end{theorem}

The geometric content of Theorem \ref{continuation} is as follows.
Let $\{\gamma_\lambda(z)\}_{z,\lambda}$ be a continuous family of
paths contained in the leaves of  $\mathcal{F}_\lambda$, and
connecting $z\in\sigma$ to $\tau$. 

For $z\in\sigma\cap
\mathbb{R^+_*}$ we suppose that $\gamma_\lambda(z)$ is the real orbit of $\mathcal{F}_\lambda$
contained in the first quadrant $x\geq 0, y\geq 0$, and connecting $z$ to $\tau$,
see fig.\ref{fig1a} (ii). The above Theorem claims that this family of orbits allows an extension to a  continuous family of
paths $\{\gamma_\lambda(z)\}_{z,\lambda}$, contained in the leaves of  $\mathcal{F}_\lambda$, and
connecting $z\in\sigma_\bullet$ to $\tau_\bullet$. The family is defined for all $(\lambda, \rho, \varphi)$ which belong to the domain (\ref{domain}). Each path starts at $z$ and terminates at a unique point on $\sigma$, denoted $\mathcal{D}_\lambda(z)$. 
Although the
paths 
$\{\gamma_\lambda(z)\}_{z,\lambda}$ are not unique, their 
relative homotopy classes
are uniquely defined. 

\subsection{Monodromy of the Dulac map and holonomy of separatrices}
To the axes $\{x=0\}, \{y=0\}$ parameterized by $y$ and $x$, we associate holonomy maps
$$
h^\lambda_\sigma : \sigma \rightarrow \sigma,\quad h^\lambda_\tau : \tau \rightarrow \tau
$$
defined by two closed paths contained in the axes $ \{x=0\}$ and $ \{ y=0 \}$ and based at $(0,c_1)$, $(c_2,0)$ respectively.
We shall make the convention, that each closed path   makes one turn around the origin of the  axe in which it is contained, in a positive direction (recall that the axes are parameterized by $y$ and $x$ respectively). It is easily seen that in the case of a linear foliation of the form
\begin{equation}
\label{linear}
x dy + \alpha \,y dx = 0, \quad \alpha \in \R^+
\end{equation}
we have
\begin{equation}
\label{linearmap}
\mathcal{D}_\alpha : x \mapsto  y= c_1c_2^{-\alpha} x^\alpha, \quad h_\sigma: x \mapsto x e^{-2\pi i /\alpha}, \quad
h_\tau: y \mapsto y e^{-2\pi i \alpha} .
\end{equation}
In the general case of a nonlinear foliation of the form (\ref{foliation}) the Dulac map $\mathcal{D}_\lambda$ is only asymptotic to $c_1c_2^{-\alpha} x^\alpha$, while the holonomy maps are analytic in $x,y,\lambda$ and
\begin{equation}
\label{nonlinear}
 \quad h_\sigma^\lambda: x \mapsto x e^{-2\pi i /\alpha}+ \dots, \quad
h_\tau^\lambda : y \mapsto y e^{-2\pi i \alpha}+ \dots, \alpha = \alpha(\lambda) .
\end{equation}
The Dulac map $\mathcal{D}_\lambda$ is a transcendental multi-valued map.
For $x>0$ let $\mathcal{D}_\lambda (e^{2\pi i} x)$ be the result of the analytic continuation of $\mathcal{D}_\lambda$ along an arc of radius $x$  and angle $2\pi i$. Similarly, for $y>0$ let $\mathcal{D}_\lambda (e^{2\pi i} y)$ be the result of the analytic continuation of $\mathcal{D}_\lambda$ along an arc of radius $y$  and angle $2\pi i$. 
\begin{lemma}
\label{monodromy}
For every sufficiently small $x>0$, $y>0$, $ |\lambda|$ holds
$$
 h^\lambda_\tau\circ  \mathcal{D}_\lambda (e^{2\pi i} x)  =  \mathcal{D}_\lambda(x), \quad
 h^\lambda_\sigma\circ  \mathcal{D}_\lambda^{-1} (e^{2\pi i} y)  =  \mathcal{D}_\lambda^{-1}(y) .
$$
\end{lemma}
{\bf Proof.} 
Consider, instead of $\mathcal{D}_\lambda$ the underlying path $\gamma_\lambda$. 
The  loop $\gamma_\lambda(e^{2\pi i} x)$ has the same origin as $\gamma_\lambda( x)$ so they can be composed and the resulting loop 
$\tilde{\gamma}_\lambda(y)$ starts at $y= \mathcal{D}_\lambda (e^{2\pi i} x)\in \tau$ and terminates at $\mathcal{D}_\lambda ( x)\in \tau$.
In the special linear case (\ref{linear}) with $\alpha=1$ the foliation is a fibration, the paths $\gamma_\lambda(.)$ represent relative cycles in the fibers of $xy$, and the path $\tilde{\gamma}_\lambda(y)$ is  closed and represents a  vanishing cycle. The claim of Lemma \ref{monodromy}
 is then the classical Picard-Lefschetz formula. In the general case the result follows "by deformation". Indeed, in the linear case  (\ref{linear}) with $\alpha=1$ the family of closed paths $\{\tilde{\gamma}_\lambda(y)\}_y$ is defined for all sufficiently small $y$, and
 $ \tilde{\gamma}_\lambda(0) \subset \{y=0\}$ is a closed path which makes one turn around the origin in the axe $\{y=0\}$ in a positive direction. Note that the paths $\tilde{\gamma}_\lambda(0)$ are bounded away from the origin in $\C^2$. 
   It follows that  $ \tilde{\gamma}_\lambda$ defines 
 the holonomy $h^\lambda_\tau$ of the separatrix $\{y=0\}$, and this property holds true also
 for every sufficiently small deformation of (\ref{linear}). The homothety $(x,y) \rightarrow (\varepsilon x, \varepsilon y)$ transforms  (\ref{foliation}) to a small deformation of (\ref{linear}) which completes the proof of the first identity
(but see also \cite{lora05}). The second identity in Lemma \ref{monodromy} is proved in a similar way. $\Box$

\subsection{The zero locus of the imaginary part of the Dulac map}
\label{zerolocus}
Consider the universal covering
\begin{equation}
\label{universal}
\C_\bullet \stackrel{\pi}{\rightarrow} \C\setminus \{0\}
\end{equation}
and the zero locus $\mathcal{H}_\lambda \subset \C_\bullet$ of the  imaginary part of the Dulac map $\mathcal{D}_\lambda$ corresponding to the domain (\ref{domain})
\begin{equation}\label{zero}
\mathcal{H}_\lambda = \{z= (\rho,\varphi)
\in \C_\bullet : \Im \mathcal{D}_\lambda(z) = 0,\,\, 0<\rho <\rho(\varphi ), \varphi \in \R \}.
\end{equation}
 In the case of a linear foliation (\ref{linear}) the zero locus is therefore a union of half-lines:
$$
\mathcal{H}_\alpha = \{z\in \C_\bullet : Im \, z^\alpha = 0 \} = \cup_{k\in \Z} \mathcal{H}_{\alpha,k}, \quad
\mathcal{H}_{\alpha,k}
= \{ (\rho, \varphi) \in \C_\bullet: \varphi = \frac{k\pi}{\alpha} \}.
$$
To describe $\mathcal{H}_\lambda$ in the case of a general  foliation of the form (\ref{foliation}), with hyperbolic ratio $\alpha(\lambda)>0$, consider the germs of real analytic sets at the origin in $\R^2=\C$
\begin{equation}
\label{clk}
C_{\lambda,k}=  \{ z\in \C=\R^2 : (h_\sigma^\lambda)^k (z)= \bar{z}   \}
\end{equation}
where $h_\sigma^\lambda$ is the holonomy map associated to the separatrix $\{x=0\}$.

\begin{figure}[htbp]
\begin{center}
 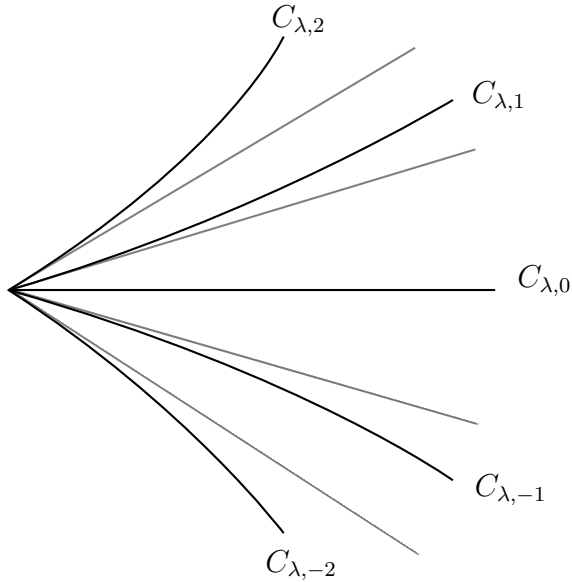
\caption{The zero locus $\mathcal{H}_\lambda$ of the imaginary
part of the Dulac map, projected on the complex plane $\C$.}
\label{dessin1}
\end{center}
\end{figure}

\begin{lemma}
\label{mainlemma}
The zero locus $\mathcal{H}_\lambda\subset \C_\bullet$ of the imaginary
part of the Dulac map in the domain (\ref{domain}) is a union of connected components $\mathcal{H}_{\lambda,k}$, indexed by $k\in \Z$.
\begin{itemize}
\item 
Each  set $C_{\lambda,k}$, (\ref{clk}), is a germ of a real analytic curve of $\R^2$, which is smooth at the origin and tangent to the line 
\begin{equation}
\label{line}
\{ z= s e^{i k\pi/\alpha(\lambda)} : s\in (\R,0)   \}
\end{equation}
 there. 
\item
Each connected component $\mathcal{H}_{\lambda,k}$ is projected on the plane $\C=\R^2$ under the map $\pi$ (\ref{universal}) to the connected component of $C_{\lambda,k}\setminus \{0\}$  tangent to the half-line  (\ref{line}), $s>0$,
at the origin.
\end{itemize}
\end{lemma}
\begin{remark}
For a general bi-holomorphic map $h_\sigma^\lambda$, vanishing at the origin,  the set $(\ref{clk})$ coincides with the origin itself. The above Lemma shows, however, that for the monodromy map $h_\sigma^\lambda$ of a saddle point of a real-analytic plane vector field, the set $C_{\lambda,k}$, (\ref{clk}), is a germ of a real analytic curve of $\R^2$, which is smooth at the origin.
The position of the connected components of $C_{\lambda,k}\setminus \{0\}$ tangent to  the half-lines 
  (\ref{line}), $s>0$ (\ref{line}), $s>0$ is
shown on fig. \ref{dessin1}.
\end{remark}


The above Lemma is the main technical result of the present paper. The
analyticity of the zero locus $\mathcal{H}_\lambda$ is responsible for the
algebraic-like behavior of the Dulac map. \\
{\bf Proof of Lemma \ref{mainlemma}}
Let $x\in \sigma\cap \R^+$ and suppose that for some $\varphi>0$, $\mathcal{D}_\lambda(e^{i\varphi}x) \in \R$. As the Dulac map is real along $\sigma \cap \R^+$, then   $\mathcal{D}_\lambda(e^{-i\varphi}x) $ is complex conjugate to $\mathcal{D}_\lambda(e^{i\varphi}x) $
 and hence
$$\mathcal{D}_\lambda(e^{-i\varphi}x)= \mathcal{D}_\lambda(e^{i\varphi}x).$$
If the point $e^{-i\varphi}x$ is seen as the inverse image of $\mathcal{D}_\lambda(e^{-i\varphi}x)$ with respect to the Dulac map $\mathcal{D}_\lambda^{-1}$, then the point $e^{i\varphi}x$ is the result of the analytic continuation of the map $\mathcal{D}_\lambda^{-1}$ along a suitable closed path of $\tau$, starting and terminating at 
$ \mathcal{D}_\lambda(e^{-i\varphi}x)$.  If we put 
$$y= \mathcal{D}_\lambda(e^{-i\varphi}x), \quad e^{-i\varphi}x= \mathcal{D}_\lambda^{-1}(y)$$
 then $e^{\pm i\varphi}x$ are two values of the multivalued map $\mathcal{D}_\lambda^{-1}(y)$ and hence, by Lemma \ref{monodromy}, they differ by a power of the monodromy $h_\sigma^\lambda$
 $$(h_\sigma^\lambda)^k( e^{i\varphi}x)=  e^{-i\varphi}x$$
 or equivalently
$$(h_\sigma^\lambda)^k( z)=  \bar{z},\quad  z=e^{i\varphi}x, \mbox{  for some  }  k\in \Z . $$ 
Clearly every such relation will correspond to a connected component $\mathcal{H}_{\lambda,k}$ of  $\mathcal{H}_\lambda$.
As $\mathcal{H}_{\lambda,k}$ is an analytic set of real dimension one, then
 $C_{\lambda,k}$ is an analytic set of dimension one too. It can be defined therefore by each of the following equivalent relations
$$
C_{\lambda,k} \subset \{ z\in \C=\R^2 : \Re [(h_\sigma^\lambda)^k (z)]=  \Re( \bar{z}) \}$$
or
$$
C_{\lambda,k} \subset \{ z\in \C=\R^2 : \Im [(h_\sigma^\lambda)^k (z)]= \Im (\bar{z})  \}.
$$
As $\frac{\partial}{\partial \bar{z}} [(h_\sigma^\lambda)^k( z)- \bar{z}] = -1$, then the linear part of the complex-analytic function
\begin{eqnarray*}
\R^2 & \rightarrow & \C \\
(z, \bar{z}) & \mapsto & (h_\sigma^\lambda)^k( z)- \bar{z}
\end{eqnarray*}
 can not be identically zero, and therefore $C_{\lambda,k}\subset \R^2$ is a real analytic curve, smooth at the origin.
It follows from (\ref{nonlinear}) that the projection of $\mathcal{H}_{\lambda,k}$ under $\pi$ on the plane $\C=\R^2$ is tangent to the half-line  (\ref{line}), $s>0$ at the origin.
$\Box$

\subsection{The argument principle}
Let $\mathbf{ D} \subset \C$ be a  relatively compact domain, with piece-wise smooth boundary, 
and $\psi : \mathbf{ D} \rightarrow \C$ an analytic function which allows a continuation to the closure $\overline{\mathbf{ D} }$. Denote by
$Z_{\mathbf{ D}}(\psi)$ the number of the zeros of $\psi$ in $ \mathbf{ D}$, counted with multiplicity. If we assume that $\psi$ does not vanish along  the border $\mathbf{ \partial D}$, then the increment of the argument $Var_{\mathbf{ \partial D}}(\arg(\psi))$ of $\psi$ along $\partial \mathbf{ D}$ oriented counter-clockwise is well defined.
 $Var_{\mathbf{ \partial D}}(\arg(\psi))$ equals the winding number of the curve $\psi(\partial \mathbf{ D}) \subset \C$ about the origin and the classical argument principle states that 
\begin{equation}
\label{argument1}
 2\pi Z_{\mathbf{ D}}(\psi) = Var_{\mathbf{ \partial D}}(\arg(\psi)) .
\end{equation}
More generally, if $\psi$ has zeros on $\mathbf{ \partial D}$, isolated or not, the variation of the argument $Var_{\mathbf{ \partial D}}(\arg(\psi))$ might be not well defined.
\begin{definition}
We say that $z\in \mathbf{ \partial D}$ is a regular zero of $\psi$ if $\psi(z)=0$, and $\psi$ allows an analytic continuation in a neighborhood of $z$ in $\C$. 
\end{definition}
If we assume that $\psi$ has only regular zeros in $\overline{\mathbf{ D} }$, then $Var_{\mathbf{ \partial D}}(\arg(\psi))$ is well defined as a sum of the increments of the argument of $\psi|_{\mathbf{ \partial D}}$ between consecutive zeros of $\psi$. Indeed, the increments are finite, because the border $\partial {D}$ is piece-wise smooth.
The argument principle can be reformulated as follows
\begin{proposition}
\label{argument}
Let $\mathbf{ D}\subset \C$ be a relatively compact domain with piece-wise smooth boundary. If $\psi : \overline{\mathbf{ D} } \rightarrow \C$ is a continuous function, analytic in $\mathbf{ D} $, and having only regular zeros in $\overline{\mathbf{ D} }$, then
\begin{equation}
\label{argument2}
2\pi Z_{\mathbf{ D}}(\psi) \leq Var_{\mathbf{ \partial D}}(\arg(\psi)) \leq 2\pi Z_{\mathbf{ D}}(\psi) + 2\pi Z_{\mathbf{\partial  D}}(\psi)
\end{equation}
\end{proposition}
{\bf Proof.}
There always exists a polynomial $P$, such that $\psi/P$ has no zeros in $\overline{\mathbf{ D} }$, so we need to verify (\ref{argument2}) for polynomials only. The set $\mathbf{ D} $ is open, connected and oriented, with piece-wise smooth boundary, which therefore has no self-intersections and has an induced orientation. The inequality
$$
0 \leq Var_{\mathbf{ \partial D}}(\arg(z)) \leq 2 \pi 
$$
allows to  "remove" the zeros along $\mathbf{ \partial D}$ and hence formula (\ref{argument1}) implies  (\ref{argument2}). $\Box$

In the present paper the first inequality in (\ref{argument2}) will be used to bound the number of the zeros $Z_{\mathbf{ D}}(.)$. For this we shall need estimates on the variation of the argument $Var_{\mathbf{ l}}(\arg(.)) $ along any  compact segment $l$ of a curve. More precisely, let $l\subset \R^2=\C$ be a compact segment of a smooth real analytic curve. Let $U\subset\C$ be an open set containing $l$ and $\psi_\lambda(z)$, $\lambda\in (\C^N,0)$, be a germ of a family of complex-analytic functions in $U$ at $\lambda=0$. For every fixed $\lambda$ such that the function $\psi_\lambda$ is not identically zero,  the variation of its argument 
$$|Var_{l}(\arg(\psi_\lambda)|$$
is well defined. 
\begin{theorem}
\label{variation} Let $l$ be a compact segment of a real analytic curve and  let
$\{\psi_\lambda\}_\lambda$ be   a family of functions analytic in a neighborhood of $l$, and depending analytically in $\lambda$.
There exists $\varepsilon_0>0$, such that 
$$
\sup_{|\lambda|<\varepsilon_0, \psi_\lambda\neq 0} |Var_{l}(\arg(\psi_\lambda)| < \infty .
$$
\end{theorem}
The above result follows from the following theorem due to Gabrielov \cite{lotozu97, gabr68}
\begin{theorem}
Let $M,N$ be real analytic varieties and consider the canonical projection
$\pi : M\times N \rightarrow N$. For every
  relatively compact semianalytic set $E\subset M\times N$,   the number of the connected components of the pre-images $\pi^{-1}(n)$ is bounded from above uniformly over $n\in N$.
\end{theorem}
{\bf Proof of Theorem \ref{variation} .} 
The number of the isolated zeros of $\psi_\lambda$ along $l$ counted with multiplicity is uniformly bounded in $\lambda$ at $\lambda=0$ (Fran\c{c}oise-Yomdin Theorem \cite{lotozu97}). On an interval between two zeros of  $\psi_\lambda(.)$ the variation of the argument divided by $2\pi$ is bounded by the number of the zeros of the imaginary part of  $\psi_\lambda$ divided by two, plus the sum of the multiplicities of the zeros of $\psi_\lambda$ at the end of the interval.
The imaginary part of  $\psi_\lambda$
 is a real analytic function in $U\subset \R^2$ and the Gabrielov Theorem implies that the number of the connected components of $\{Im( \psi_\lambda)=0\} \cap l $ 
is uniformly bounded in $\lambda$ at $\lambda=0$.$\Box$

\section{Cyclicity of one-saddle cycles}
\label{onesaddle}
\begin{figure}[htbp]
\begin{center}
 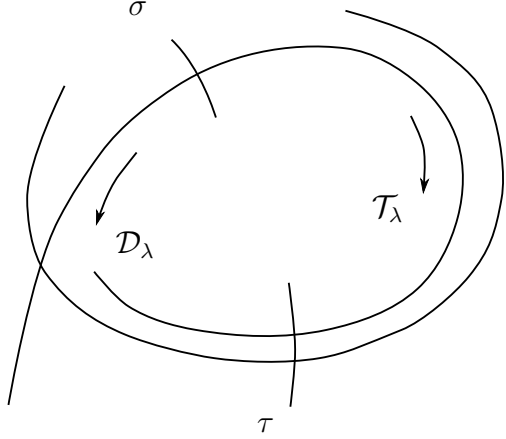
\caption{The Dulac map $\mathcal{D}_\lambda(z)$ and the transport map $\mathcal{T}_\lambda(z)$.}
\label{dessin2}
\end{center}
\end{figure}
\begin{figure}[htbp]
\begin{center}
 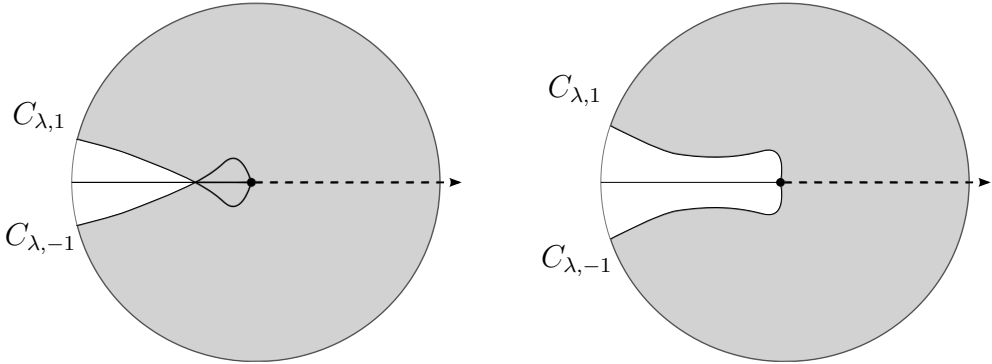
\caption{Examples of domains $\mathbf{ D}_R \subset \C_\bullet$, projected on the complex plane $\C$ under $\pi$ (\ref{universal}).}
\label{dessin4}
\end{center}
\end{figure}

Let 
$X_\lambda$, $\lambda \in(\mathbb{R}^N,0) $ be a germ of an
analytic family of analytic plane vector fields, such that $X_0$ has a one-saddle cycle (homoclinic saddle loop) $\Gamma_1$. The  first-return map associated to $\Gamma_1$ is a composition of a Dulac map $\mathcal{D}_\lambda(z) : \sigma \rightarrow \tau$ and a transport map $\mathcal{T}_\lambda(z)$, see fig. \ref{dessin2}.
We assume that the Dulac map is in a normal form as in section \ref{scontinuation}.
The limit cycles of $X_\lambda$ near $\Gamma_1$ correspond to the zeros of the displacement map
$$
 \psi_\lambda(z) =\mathcal{D}_\lambda(z) - \mathcal{T}_\lambda(z)
$$
near $z=0$.  An appropriate choice of the local coordinates on the croos-sections $\sigma$ and $\tau$ brings the transport map to the form $\mathcal{T}_\lambda(z)\equiv z$. Alternatively, we could choose simply $\sigma=\tau$ (without supposing that the Dulac map is in the normal form of section \ref{scontinuation}).
We shall bound the number of the zeros of $\psi_\lambda$ in the domain $\mathbf{ D}_R\subset \C_\bullet$ delimited by the circle
 $\{ \rho =R\} $, and the connected components $\mathcal{H}_{\lambda,1}$ and $\mathcal{H}_{\lambda,-1}$ of the zero locus of the imaginary part of the Dulac map, as it is shown on fig. \ref{dessin4}.
We shall suppose that $R>0$ is so small, that $\psi_\lambda(.)$ is analytic in $\mathbf{ D}_R$ for all $\lambda\in\R^N$,  such that $|\lambda|\leq \varepsilon_0$ (Theorem \ref{continuation}) and it is analytic even on the closure of $\mathbf{ D}_R$ except of course at $z=0$, where  $\psi_\lambda(.)$ is only continuous. Indeed, 
$$\lim_{z\rightarrow 0, z\in \mathbf{ D}_R}\mathcal{D}_\lambda(z)  =0$$
while $\mathcal{T}_\lambda(z)$ is holomorphic at $z=0$, so 
$$\lim_{z\rightarrow 0, z\in \mathbf{ D}_R}\psi_\lambda(z)  =c(\lambda)$$
where $c(\lambda)$ is analytic and $c(0)=0$. 
If the family of functions $\psi_\lambda$ is sufficiently general, then $c(\lambda)\not\equiv 0$, and in the case when $c(\lambda)\equiv 0$ we may replace $\psi_\lambda$ by the new family $\psi_\lambda+ \lambda_{N+1}$, $\lambda_{N+1}\in \R$, for which the limit at $z=0$ is the parameter $\lambda_{N+1}$. After this preparation, we may prove the finite cyclicity of the homoclinic loop $\Gamma_1$. For this we apply Proposition \ref{argument}  (the argument principle) to the family of functions $\psi_\lambda$ in the domain $\mathbf{ D}_R$. In the course of the computation, it will be supposed that $R>0$ is sufficiently small, $\varepsilon_0$ is sufficiently small with respect to $R$, and $\lambda$ is such that $|\lambda|< \varepsilon_0$. We may encode this choice of the parameters by the "physical" notation
\begin{equation}
\label{constants}
0< |\lambda|< \varepsilon_0<< R << 1 .
\end{equation}

The hyperbolic ratio of the saddle point will be not bigger than one only in a suitable semi-analytic set in the parameter space, and will be bigger than one in another (complementary) semi-analytic set. After eventual exchanging  of $\sigma$ and $\tau$, it will be also supposed that the hyperbolic ratio of the saddle point is not bigger than one for all parameter values.

Along the circle $\{z: |z|=R \}$  with angle close or strictly less than $2\pi$
the variation of the argument of $\psi_\lambda$ is uniformly bounded in $\lambda$ (Theorem \ref{variation}). 

Along the curve $C_{\lambda,1}$
the imaginary part of $\psi_\lambda$ equals the imaginary part of the 
transport map $- \mathcal{T}_\lambda(z)=-z$. Therefore the zeros of $ \Im (\psi_\lambda)$ along $C_{\lambda,1}$ are exactly
 the intersection points of $C_{\lambda,1}$ and the segment $(-R,0)$. 
According to 
Lemma \ref{mainlemma}  we have
\begin{equation}
\label{1loop}
C_{\lambda,1}\cap \R = \{ x\in \R: h_\sigma^\lambda(x)=x\} = C_{\lambda,-1}\cap \R .
\end{equation}
As $h_\sigma^\lambda(x)$ is an analytic family of analytic functions, then by Gabrielov's theorem, the number of such fixed points is uniformly bounded in $\lambda$ on $[-R,0]$. To conclude, we have only to check that the family $\{\psi_\lambda\}_\lambda$ has regular zeros along the border of the domain $\mathbf{ D}_R$. This is indeed the case, when $c(\lambda)\neq0$, as $\psi_\lambda(0)=c(\lambda)$. We conclude that the number of isolated zeros of the family of functions
$$
\{\psi_\lambda: c(\lambda) \neq 0, |\lambda|\leq \varepsilon_0\}
$$
in the domain $\mathbf{ D}_R$ is uniformly bounded by some integer, say $C$. Finally, note that the condition $c(\lambda) \neq 0$ can be removed. Indeed, if for some $\lambda_0$, $|\lambda_0| \leq \varepsilon_0$, $c(\lambda)=0$, the function $\psi_{\lambda_0}$ has at least $C+1$ zeros in  $\mathbf{ D}_R$, then it has at least $C+1$ zeros in $\mathbf{ D}_R$ in a sufficiently small neighborhood of $\lambda_0$, in contradiction with the preceding estimate.

To resume, we proved the following classical result\\
{\bf Theorem (Roussarie\cite{rous86, rous89, rous98a})}
\emph{Every homoclinic saddle loop (a one-saddle cycle) occurring in an analytic finite-parameter family of plane analytic vector fields, may generate
 no more than a finite number of limit cycles within the family.}

Let us note that our method, exactly as the Roussarie's Theorem  allows to compute more precisely the cyclicity of $\Gamma_1$. 
We shall not enter into details here. Just to illustrate this, note that if the hyperbolic ratio $\alpha(0)$ is strictly bigger than one, then  the overall increase of the argument of the displacement map along the border of $\mathbf{ D}_R$ is strictly less than $2\pi$ (this computation is omitted) and the cyclicity of $\Gamma_1$ is zero.

\section{The Petrov trick}
\label{petrovtrick}
The content of this section is not necessary for the proof of our main result Theorem \ref{mainth}, but it aims to shed some light on the origin of the method, used to bound the limit cycles near the saddle loop in the preceding section.

With the same notations as in section \ref{onesaddle}, consider the analytic family of analytic vector fields $$X_\lambda, \;\lambda= (\lambda_1, \dots, \lambda_N) \in(\mathbb{R}^N,0) $$ defining a  holomorphic foliation $\mathcal{F}_\lambda$ of the form
$$
\mathcal{F}_\lambda = \{ dH+ \lambda_1\omega_\lambda=0\}, \omega_0 \neq 0
$$
 where $H$ is a function and
$\omega_\lambda$ is an analytic family of differential one-forms,  both analytic in a neighborhood of the saddle loop $\Gamma_1$.
For definiteness, we put the saddle point at the origin in $\R^2$, so $dH(0)=0$. We shall further suppose that the saddle loop $\Gamma_1$ is contained in the level set $\{H(x,y) = 0 \}$, and the interior of $\Gamma_1$ is filled up by a continuous family of periodic orbits $\gamma_0(h)\subset  \{H(x,y) = h \}$, parameterized by $h> 0$, where $h=H(x,y)|_\sigma$ is the restriction of $H$ on the cross-section $\sigma$.
The displacement map  is approximated by the usual Poincar\'e-Pontryagin formula as follows
\begin{equation}
\psi_\lambda(h)= \lambda_1 \int_{\gamma_0(h)} \omega_\lambda + o(\lambda_1), 
\end{equation}
where $o(\lambda_1)/\lambda_1$ tends to zero as $\lambda$ tends to zero, uniformly in  $h $ in every compact interval in which the displacement map is defined. The zeros of $\psi_\lambda(.)$ correspond to limit cycles and, at least far from $\Gamma_1\subset \{H(x,y) = 0 \}$, they are approximated there by the zeros of the complete Abelian integral 
$$
h \mapsto I_\lambda(h)=\int_{\gamma_0(h)} \omega_\lambda, \;\; h \geq 0 .
$$
We make the assumption (actually justified by the Roussarie's theorem \cite{rous86}), that this is so also in a neighborhood of $h~=~0$ (corresponding to limit cycles close to the saddle loop $\Gamma_1$). Thus, it makes a sense  to prove the finiteness of the maximal number of the zeros of the Abelian integral $I_\lambda(h)$, which tend to 
$h=0$ as $\lambda $ tends to the origin in the parameter space. 
This follows of course from a well known  general result of Varchenko and Khovansky. We shall use, however, a different idea due to G.S. Petrov \cite{petr88b},  who showed that the analogous global problem for complete elliptic integrals of second kind is of algebraic nature. This observation has been used in several papers by Petrov to evaluate the precise number of zeros of complete elliptic integrals, and hence of limit cycles of perturbed Hamiltonian vector fields, see for instance {\.Z}o{\l}adek \cite[section 6]{zola06}. 
We are ready to describe the local version of the Petrov method. 
\begin{quotation}
Consider the sector $$S_R= \{z= \rho e^{i\varphi} \in \C: 0< \rho< R,\; 0< \varphi < 2\pi \} .$$ For a fixed sufficiently small $R>0$ and all sufficiently small $\|\lambda\|$ the Abelian integral $I_\lambda(z)$ allows an analytic continuation in $S_R$. To bound the number of its zeros on $S_R$ (and hence on $(0,R)$) we apply the argument principle to the domain $S_R$. Along the circle $\{\rho = R \}$ the increase of the argument of $I_\lambda$ is bounded uniformly in $\lambda$ (due to Gabrielov's theorem). Along the segment $[-R,0]$ the Abelian integral allows two analytic continuations $I_\lambda^\pm(h)$. As $I_\lambda(.)$ is real-analytic on $(0,R)$ then
$$
I_\lambda^+(h) = \overline{I_\lambda^-(h)}, \; h\in (-R,0)
$$
and by the Picard-Lefschetz formula
\begin{equation}
\label{impart}
2 \sqrt{-1}\Im I_\lambda^+(h) =  I_\lambda^+(h) - I_\lambda^-(h) = \int_{\delta(h)} \omega_\lambda,  \; h\in (-R,0)
\end{equation}
where $\delta(h) \subset  \{H(x,y) = h \}$ is a continuous family of cycles, vanishing at the origin as $h$ tends to zero. 

The imaginary part of $I_\lambda(h)$ on $(-R,0)$ is therefore an analytic function, and by Gabrielov's theorem again, its zeros are uniformly bounded in $\lambda$ on the closed interval $[-R,0]$. This implies that the increase of the argument of $I_\lambda(h)$ on $(-R,0)$ is also uniformly bounded in $\lambda$ which combined to the argument principle shows the finiteness of the maximal number of zeros 

\end{quotation}
The proof of the finite cyclicity of the one-saddle loop from the preceding section, may be seen as a generalization of the Petrov method. Indeed, the Picard-Lefschetz formula corresponds to the claim of Lemma 1, and by Lemma 2 the zeros the analytic Abelian integral
(\ref{impart}) correspond to the fixed points (complex limit cycles) of the holonomy map $h_\sigma^\lambda$ of the separatrix. As it is well known, the holonomy map of a separatrix is analytic, which implies the finite cyclicity of the saddle loop $\Gamma_1$.

\section{Cyclicity of two-saddle cycles}
\label{mainsection}

The main result of the paper is the following
\begin{theorem}
\label{mainth}
Every heteroclinic saddle loop (a two-saddle cycle) occurring in an analytic finite-parameter family of plane analytic vector fields, may generate
 no more than a finite number of limit cycles within the family.
\end{theorem}

\begin{figure}
\begin{center}
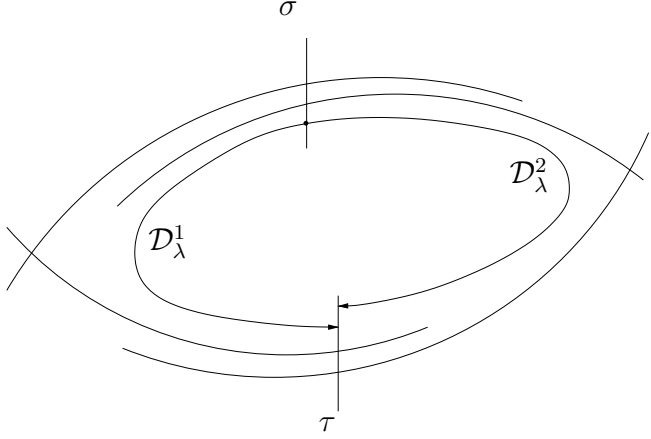
\end{center} \caption{The Dulac maps $\mathcal{D}^1_\lambda$ and
$\mathcal{D}^2_\lambda$ }
\label{fig4}
\end{figure}

\begin{figure}[htbp]
\begin{center}
 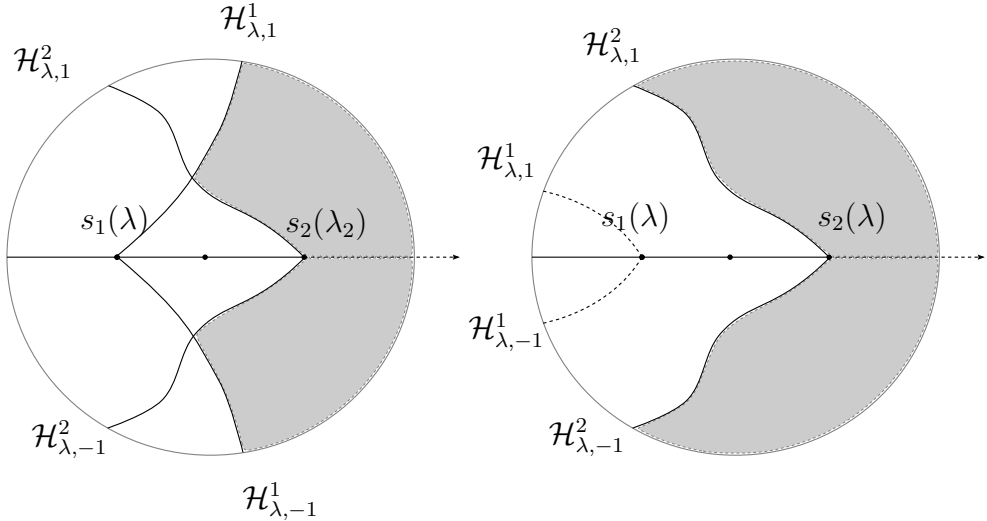
\caption{The domain $\mathbf{ D}_R \subset \C_\bullet$ projected on the complex plane $\C$ in the case $ \alpha_2(0)>1$.}
\label{dessin5}
\end{center}
\end{figure}

\begin{figure}[htbp]
\begin{center}
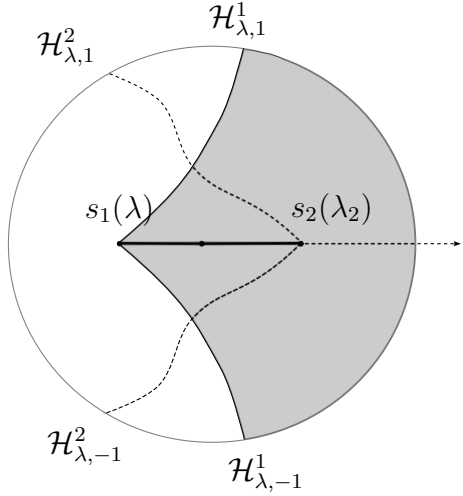
\caption{The domain $\mathbf{ D}_R \subset \C_\bullet$ projected on the complex plane $\C$ in the case $ \alpha_2(0)\leq 1$, $ \alpha_1(0)\geq 1$.}
\label{dessin6}
\end{center}
\end{figure}

Using the notations used in the preceding sections,   suppose that the vector field $X_0$ has a
two-saddle loop $\Gamma_2$. Consider the Dulac maps
 $$
\mathcal{D}^i_\lambda : \sigma \rightarrow \tau, \; i=1,2
$$
associated to the corresponding foliation, as on
fig.\ref{fig4}. 
Each map $\mathcal{D}^i_\lambda$ is a composition of a "local" Dulac
map (as in section \ref{dulacsection}) and two real-analytic transport maps. From
this it follows that Lemma \ref{mainlemma} applies to  $\mathcal{D}^i_\lambda$,
$i=1,2$, too. 
From now on we choose a real-analytic local variable $z$ on the cross-section $\sigma$ thus identifying $\sigma$ to an open disc centered at $0\in \C$. We shall also suppose that $0=\sigma\cap\Gamma_2$.
 The functions $\mathcal{D}_\lambda^i(z)$, $i=1,2$ are multivalued on the cross-section $\sigma$ and have 
critical points at $s_i(\lambda)\in \mathbb{R}$, $s_i(0)=0$, respectively. The functions
$s_i$ are real-analytic.  
The limit cycles of $X_\lambda$ near $\Gamma_2$ correspond to the zeros of the displacement map
$$
 \psi_\lambda(z) =\mathcal{D}_\lambda^1(z) - \mathcal{D}_\lambda^2(z)
$$
near $z=0$.  Let $\alpha_i(\lambda)>0$, $i=1,2$ be the hyperbolic ratios of the saddles.  We shall suppose, upon exchanging eventually the roles of $\sigma$ and $\tau$,  that $\alpha_1(0) \alpha_2(0) \geq 1$. Denote the zero loci of the imaginary parts of 
the Dulac maps $\mathcal{D}_\lambda^1(z)$,  $\mathcal{D}_\lambda^2(z)$ by
$\mathcal{H}_{\lambda}^1$ and $\mathcal{H}_{\lambda}^2$ respectively. 
We shall bound the number of the zeros of $\psi_\lambda$ in the complex domain
$\mathbf{D}_R$ of the universal covering of $\C\setminus \{s_1(\lambda), s_2(\lambda) \}$ defined as follows (without loss of generality we assume that $s_1(\lambda) \leq s_2(\lambda)$). 
\begin{itemize}
\item 
if $ \alpha_2(0)>1$, the domain $\mathbf{D}_R$ is  bounded
by the circle
\begin{equation}\label{sr}
S_R= \{z: |z| = R \},
\end{equation}
and by
$$\mathcal{H}_{\lambda,1}^1, \mathcal{H}_{\lambda, -1}^1, \mathcal{H}_{\lambda,1}^2, \mathcal{H}_{\lambda,-1}^2$$
as it is shown on fig.\ref{dessin5}.
\item
if
$\alpha_2(0)\leq 1$ then necessarily $\alpha_1(0) \geq 1$. The domain $\mathbf{D}_R$ is  bounded
by the circle
$S_R$,
by the interval $[s_1(\lambda),s_2(\lambda)]$, and by $\mathcal{H}_{\lambda,1}^1, \mathcal{H}_{\lambda, -1}^1$, as it is shown on fig.\ref{dessin6}.
\end{itemize}

In the course of the proof the parameters $R$ and $\lambda$ will be chosen as in the one-saddle case: the constant $R$ will be sufficiently small, $\varepsilon_0>0$ will be sufficiently small with respect to $R$, and 
$\lambda \in \R^N$ will be such that
$|\lambda|<\varepsilon_0$, see (\ref{constants}). Like in section \ref{onesaddle} we shall suppose, without loss of generality, that the analytic functions $c_1(\lambda), c_2(\lambda)$ where
$$
\lim_{z\rightarrow s_1(\lambda), z\in \mathbf{ D}_R}\psi_\lambda(z)  =c_1(\lambda),
\lim_{z\rightarrow s_2(\lambda), z\in \mathbf{ D}_R}\psi_\lambda(z)  =c_2(\lambda)
$$
are not identically zero. This will guarantee that for generic values of $\lambda$ the displacement map will have only regular zeros in the closure of $\mathbf{D}_R$, so the argument principle  (Proposition \ref{argument})  can be applied.

{\bf Proof of Theorem \ref{mainth}.}
It follows from the definition of the domain $\mathbf{ D}_R \subset \C_\bullet$  that  the displacement map 
$\psi_\lambda(z)$
is analytic there. To count the zeros (corresponding to real and complex  limit cycles) of the displacement map in 
$\mathbf{ D}_R $ we apply Proposition 1 (the argument principle) to the family of functions $\psi_\lambda$.
To evaluate the  variation of the argument of the displacement map  along the border of $\mathbf{D}_R$ we 
repeat the arguments of section \ref{onesaddle}. 

Consider first the case $ \alpha_2(0)>1$,  fig.\ref{dessin5}.  The connected component of the zero locus of the imaginary part of $\mathcal{D}_\lambda^2$ which is tangent to the line $\varphi=  \pi /\alpha_2(\lambda)$ through $s_2(\lambda)$ intersects the circle $S_R$ transversally, and along this circle  the variation of the argument of $\psi_\lambda$ is uniformly bounded in $\lambda$ (Theorem \ref{variation}). 
 The imaginary part of $\psi_\lambda(z)$ restricted to  $\mathcal{H}_{\lambda}^1$ equals the imaginary part of $-\mathcal{D}_\lambda^2$ and hence $\Im \psi_\lambda$ vanishes along $\mathcal{H}_{\lambda,1}^1, \mathcal{H}_{\lambda,-1}^1$ exactly at the intersection points
 $$
 \mathcal{H}_{\lambda,1}^1 \cap  \mathcal{H}_{\lambda,1}^2, \quad  \mathcal{H}_{\lambda,-1}^1 \cap  \mathcal{H}_{\lambda,-1}^2 .
 $$
According to Lemma \ref{mainlemma} these intersection points are the solutions of the equation
\begin{equation}
\label{8loop}
h^\lambda_2(z) = h^\lambda_1(z)
\end{equation}
where $h^\lambda_1, h^\lambda_2$ are the holonomies of the separatrices intersecting $\sigma$ and related to the saddle points $s_1(\lambda)$ and $s_2(\lambda)$. By  Gabrielov's theorem, the number of such fixed points is uniformly bounded in the disc $\{z: |z| < R \}$. 

 Consider now the second case $\alpha_2(0)\leq 1$, $\alpha_1(0) \geq 1$, see fig.\ref{dessin6}. Along this circle  $S_R$ the variation of the argument of $\psi_\lambda$ is uniformly bounded in $\lambda$ (Theorem \ref{variation}). Along the interval $[s_1(\lambda),s_2(\lambda)]$ the imaginary part of $\mathcal{D}_\lambda^1$ vanishes identically, and the imaginary part of $\psi_\lambda(z)$ restricted to  this interval equals the imaginary part of $-\mathcal{D}_\lambda^2$. 
 Therefore the zeros of $ \Im (\psi_\lambda)$ along $[s_1(\lambda),s_2(\lambda)]$ are exactly
 the intersection points of $\mathcal{H}_{\lambda,1}^2$ and $[s_1(\lambda),s_2(\lambda)]$. By Lemma \ref{mainlemma}, and  like in (\ref{1loop}), these intersection points are the solution of the equation $$h^\lambda_2(z) =z$$ where $h''$ is the holonomy of the separatrix intersecting $\sigma$ and related to the saddle points $s_2(\lambda)$. By  Gabrielov's theorem, the number of such fixed points is uniformly bounded. Finally,  the zeros of $ \Im (\psi_\lambda)$ along 
 $\mathcal{H}_{\lambda,1}^1$ and $\mathcal{H}_{\lambda,-1}^1$ are evaluated as in the case $ \alpha_2(0)>1$.
  This completes the proof of Theorem \ref{mainth}.$\Box$

\section{Concluding remarks.} 

\begin{figure}[htbp]
\begin{center}
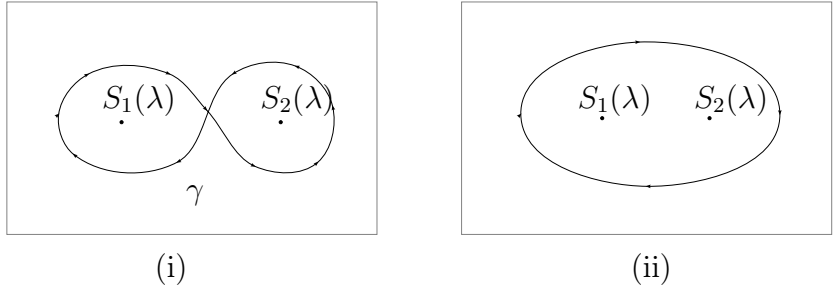
\caption{The figure eight-loop $\gamma$.}
\label{dessin8}
\end{center}
\end{figure}

The identity (\ref{8loop}) which determines complex limit cycles "responsible" for the cyclicity of the double loop $\Gamma_2$ is the main new ingredient in the proof with respect to the one-saddle case. Indeed, solutions of (\ref{8loop})  are fixed points of the holonomy $h^\lambda_2\circ (h^\lambda_1)^{-1}$ which, for $\lambda=0$, is generated by a closed loop $\gamma$ contained in the complexified separatrix of $\Gamma_2$ intersecting the cross-section $\sigma$. The topological type of this separatrix near $\Gamma_2$ is a disc with two punctures, corresponding to the two saddle points $S_1(\lambda)$ and $S_2(\lambda)$. Clearly $\gamma$ makes one turn around each of them, but depending on the orientation we have two possibilities shown on fig. \ref{dessin8} (i) and (ii). A simple computation on a model example shows that the loop $\gamma$ associated to the holonomy $h^\lambda_2\circ (h^\lambda_1)^{-1}$ is the figure eight-loop on fig. \ref{dessin8} (i). The reader will recognize in the loop $\gamma$ a key ingredient in the proof of the local boundedness of the number of  zeros of pseudo-Abelian integrals in \cite{bmn09,bmn11}.

Although the result of Theorem \ref{mainth}  is existential, the proof we use leads to effective upper bounds on the number of the bifurcating limit cycles. This possibility is explored in 
 \cite{gavr11}, where we show that the cyclicity of a Hamiltonian two-loop is bounded by the number of the zeros of \emph{a pair} of associated Abelian integrals, a phenomenon which also explains the appearance of alien limit cycles in \cite{duro06}. 
  \begin{figure}[htbp]
\begin{center}
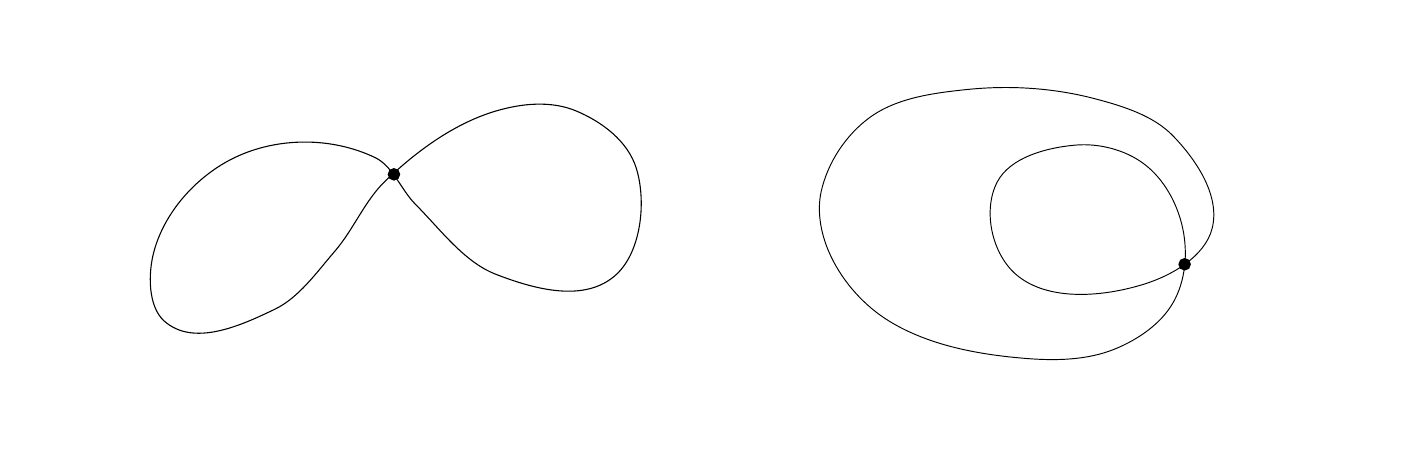
\caption{Hyperbolic planar polycycles with finite cyclicity.}
\label{dessin7}
\end{center}
\end{figure}

It worth noting, that our finitness result holds true, with the same proof,  for  other hyperbolic polycycles (on the plane or on an analytic surface), as those shown on fig. \ref{dessin7}.  

\paragraph{Acknowledgments.}
The author thanks Marcin Bobie\'nski for the stimulating discussions, as well to the anonymous referees for the valuable suggestions. 

\newpage

\def\cprime{$'$} \def\cprime{$'$} \def\cprime{$'$} \def\cprime{$'$}
  \def\cprime{$'$} \def\cprime{$'$} \def\cprime{$'$}

\end{document}